\documentclass[12pt]{elsarticle}

\usepackage[top=99pt, bottom=75pt, left=72pt, right=70pt]{geometry}
\usepackage[utf8]{inputenc}
\usepackage{amsthm}
\usepackage{amsmath}
\usepackage{amssymb}
\usepackage{amscd,amsxtra}
\usepackage{MnSymbol}
\usepackage{mathtools}
\usepackage{tensor}
\usepackage{faktor}
\usepackage{dsfont}
\usepackage[all]{xy}
\usepackage{graphicx}
\usepackage[rgb]{xcolor}
\usepackage{tikz, tikz-cd}
\usetikzlibrary{arrows, matrix, intersections, math, calc, decorations.pathmorphing, decorations.markings, positioning}
\usepackage{enumitem}
\usepackage{hyperref}
\definecolor{R}{RGB}{255, 0, 34}
\definecolor{B}{RGB}{0, 85, 238}
\hypersetup{colorlinks = true,
            linkcolor = R,
            urlcolor  = B,
            citecolor = B,
            anchorcolor = B}

\newtheorem{theorem}{Theorem}

\numberwithin{theorem}{subsection}
\numberwithin{proposition}{subsection}
\numberwithin{lemma}{subsection}
\numberwithin{claim}{subsection}
\numberwithin{corollary}{subsection}
\numberwithin{conjecture}{subsection}
\numberwithin{definition}{subsection}
\numberwithin{remark}{subsection}
\newtheorem{example}[theorem]{Example}

\newcommand{\beq}{\begin{equation}}
\newcommand{\eeq}{\end{equation}}
\newcommand{\beqa}{\begin{eqnarray}}
\newcommand{\eeqa}{\end{eqnarray}}
\newcommand{\beaa}{\begin{eqnarray*}}
\newcommand{\ben}{\begin{eqnarray*}}
\newcommand{\eaa}{\end{eqnarray*}}
\newcommand{\een}{\end{eqnarray*}}

\makeatletter
\def\ps@pprintTitle{%
  \let\@oddhead\@empty
  \let\@evenhead\@empty
  \def\@oddfoot{\reset@font\hfil\thepage\hfil}
  \let\@evenfoot\@oddfoot
}

\journal{}

\begin{document}

\begin{frontmatter}



\title{\textbf{On Physical Mathematics: \\
an approach through Gilles Ch\^{a}telet's philosophy}} 


\author{John Alexander Cruz Morales}

\ead{jacruzmo@unal.edu.co} 
\ead{johnalexander.cruzmorales@csg.igrothendieck.org}

\affiliation{organization={Universidad Nacional de Colombia, Departamento de Matemáticas},
            addressline={Ciudad Universitaria}, 
            city={Bogotá},
           country={Colombia}}

\affiliation{organization={Istituto Grothendieck, Centro di Studi Grothendieckiani},
            addressline={Corso Statuto 24}, 
            postcode={12084},
            city={Moldov\`{i}},
            country={Italy}}                
            
\begin{abstract}
\small{Starting from Greg Moore's description about Physical Mathematics, a framework is proposed in order to understand it, based on Gilles Ch\^atelet's philosophy. It will be argued that Ch\^atelet's ideas of inverting, splitting, augmenting and virtuality are crucial in the discussion about the nature of Physical Mathematics. Along this line, it will be proposed that mirror symmetry is a natural study case to test Ch\^atelet's ideas in this context. This should be considered as a first step in a long term project aiming to study the relations among mathematics, physics and philosophy in the construction of a global understanding of the structure of the universe, as it was envisioned by Grothendieck in the late 80's of the last century and it was started to be developed independently by Ch\^atelet in the beginning of the 90's. The main suggestion of the essay is that it is in the relations between mathematics, physics and philosophy that new knowledge arises.} 
\end{abstract}




\end{frontmatter}


\begin{flushright}
\footnotesize{\textit{``Nothing is more fertile, all mathematicians know, than these obscure analogies, these murky reflections of one theory in another, these furtive caresses, these inexplicable tiffs; also nothing gives as much pleasure to the researcher. A day comes when the illusion vanishes; presentiment turns into certainty; twin theories reveal their common source before disappearing; as the \emph{Gita} teaches, one attains knowledge and indifference at the same time.'' \\
― Andr\'e Weil, De la m\'etaphysique aux math\'ematiques.}}

\end{flushright}

\section{The triangle: mathematics, physics, philosophy in a nutshell}

Andr\'e Weil's quote in the preamble also appears in the introduction of Ch\^atelet treatise  \emph{Les enjeux du mobile. Math\'ematique, physique, philosophie} \cite{chatelet1}. I have reproduced it here because it exemplifies the spirit in which this short note is written. This text will try to explore some \textbf{``common sources''} and some \textbf{``obscure analogies''} between mathematics and physics and their further interaction with philosophy.\\

More concretely, one of the goals of this work is to show how Gilles Ch\^atelet's ideas can be used to understand from a philosophical point of view the emergent field called \textbf{Physical Mathematics}, and also, based on this, to give the first steps in a long term project aiming to study the relationships among physics, mathematics and philosophy in order to produce new insights that allow us to improve our understanding of the fundamental questions about the structure of the universe. This project fits into the one proposed by Grothendieck in \cite{grothendieck1, grothendieck2}, where he proposed that questions, like \textbf{what is a space?}, must be considered in an unified way from a mathematical, physical and philosophical point of view. In fact, in Grothendieck's own words, see \cite{grothendieck1}, regarding the question on the nature of the notion of space, he writes: \\

\emph{``In summary, I expect that the awaited renewal (if one is yet to come...) will need to come from someone who is a mathematician at heart, and well-versed in physics' great challenges, rather than from a physicist. And most importantly, this person will need to have the ``philosophical openness'' to grasp the heart of the problem - for this problem is not a technical one, but rather a fundamental problem in the ``philosophy of nature'' ''}. \\

This part of Grothendieck's astonishing work is almost unknown for non-specialist and little known for specialists. As a by product, this essay pretends to show this set of Grothendieck's ideas where the traditional borders between different fields of knowledge (a.k.a mathematics, physics and philosophy) are blurred, allowing the emergency of interesting and complementary perspectives that are needed, specially nowadays, when the hyper-specialization is the new rule. Therefore, Grothendieck's proposal should be understood as a modern call of the spirit that dominated the intellectual endeavours in the past. The above mentioned current hyper-specialization, which started at the end of the 19th century and consolidated through the 20th century, produces the illusion that global attempts to understand the universe are sort of nonsense and no longer useful enterprises. I really disagree with this approach and fortunately notable examples like Grothendieck and Ch\^atelet, in recent decades, show that there is still space for this kind of dreams and visions. They are two of the giants whose shoulders I will use to stand up. \\

Of course, the actual content and intention of this note are far more modest. Motivated by the description of \emph{Physical Mathematics} due to Gregory Moore, one of the leading figures in the field, which propose a new way to understand the traditional relations between mathematics and physics, I will argue that the main features of this description lead to a context that can be better understood using Ch\^atelet's work on the \textbf{philosophy of physico-mathematics}. In fact, I would like to argue that the physico-mathematics in Ch\^atelet is a version of the Physical Mathematics in Moore and vice-versa. This justifies my choice of Ch\^atelet's philosophy as a right setting for approaching Moore's Physical Mathematics. \\

The text is organized as follows. In the first section, I will discuss Moore's characterization of Physical Mathematics and will point out its main features, namely, inversion, differentiation and virtuality. The main argument here is that these three aspects suggest the necessity of a philosophical framework that takes care of them in order to allow the construction of an understanding about what Physical Mathematics is with solid philosophical grounds\footnote{Moore, in his description on Physical Mathematics, claims to start with a philosophical approach. I think this is the right way to start. However, the question about what would be the more adequate philosophical background in order to discuss about Physical Mathematics remained open. I think this a serious question and this text is an attempt to address this important issue.}. In the second section, I will discuss Ch\^atelet's ideas about physico-mathematics and, based on the analysis done by Zalamea in \cite{zalameachat}, will argue that the three processes appearing in Moore's definition correspond to the three main gestures proposed by Ch\^atelet in his philosophy. It is important to remark that I will not pretend to explain Ch\^atelet's global system in detail. This is beyond the scope of this text. For a more detailed presentation, see \cite{zalameachat}. \\

I hope that the philosophical approach started in this paper will help to tackle the difficulties that Physical Mathematics has faced in the community of mathematicians and physicists. As Moore points out in \cite{moore1}: \\

\emph{``Physical Mathematics is sometimes viewed with suspicion by both physicists and mathematicians. On the one hand, mathematicians regard it as deficient, for lack of proper
mathematical rigour...On the other hand, the relative lack of reliance of Physical Mathematics on laboratory experiments is viewed - with some justification - as dangerous by many physicists''.}\\

From my point of view, the lack of a systematic philosophical treatment is one of the sources of such suspicions. In this direction, one needs the three vertices of the triangle (mathematics, physics, philosophy) to see the actual power of Physical Mathematics. A better picture would be to consider the vertices as providing local charts and, in order to have a global view, one needs to glue them together in a coherent way. This is reached with a fluid dialogue among mathematics, physics and philosophy in the form of a philosophy of Physical Mathematics. With this goal in mind, finally, at the end of the text, I will propose some ideas on how Physical Mathematics understood \emph{\'a la Ch\^atelet} might be used to increase our current understanding of the mirror phenomenon, not just from a philosophical point of view, but also mathematically and physically. Following Grothendieck's quote above, the idea I want to express is that deep down mirror symmetry is a question on natural philosophy.

\section{Emergence of Physical Mathematics: relations between mathematics and physics revisited}
The relation between mathematics and physics is quite old and certainly is at the heart of the development of modern science. It is difficult to determine whether Newton or Leibniz were physicists or mathematicians\footnote{Or even philosophers. In particular, Leibniz was a global thinker and it is difficult, if not impossible, to classify him in a small field.}. One might face the same problem when looking at the work of Euler, Lagrange or Poincar\'e. Even Riemann, who undoubtedly is a mathematician in the modern sense, liked to think of himself also as a physicist\footnote{Riemann also had very strong philosophical interests. Maybe, it is not a mere coincidence that two of the more important figures in transforming our understanding of the space, namely Riemann and Grothendieck, have also try to develop a philosophical picture of the world.}. The reason for this is that mathematics and physics are strongly tie to each other that trying to stablish a border between them, in many cases, is unclear. \\

However, in a certain sense, the relation has not been reciprocal. It was traditional to think that physics was only a receptacle of mathematical applications and it was mysterious why this happened. An outstanding example of this attitude could be found in the title of the famous essay by Wigner \emph{On the Unreasonable Effectiveness of Mathematics in the Physical Sciences} \cite{wigner}, and also encoded in the expression Mathematical Physics. In this sense, physics provides some problems about the real world that ``magically'' can be solved using mathematical tools. Mathematical intuition helps the development of physics but not the other way around, i.e., it was traditional to think that the physical intuition does not help to solve purely mathematical problems\footnote{This statement should be taken carefully since a serious historical revision of the works of Riemann or Poincar\'e, just to mention two remarkable examples, might lead to a different conclusion. I do not pretend to address this issue here. In any case, as a general belief, it was accepted that the only contribution of physics in mathematics was to provide interesting problems.}.  \\

In the early 70's of the last century, the relation between mathematics and physics seemed to be at its lowest possible point. In 1972 Freeman Dyson wrote his essay \emph{Missed Opportunities} \cite{dyson}, where he declares: \\

\emph{``As a working physicist, I am acutely aware of the fact that the marriage between
mathematics and physics, which was so enormously fruitful in past centuries, has recently
ended in divorce''.}\\

Fortunately, Dyson's pessimistic view was not shared by many people. As Moore points out, see \cite{moore1, moore2}, around the same time when Dyson made his declaration, the work of other people including Atiyah, Bott and Singer from the mathematical side and Coleman, Gross and Witten from the physical side, led to a reconstruction of the marriage with new and powerful insights. From the mathematical side, mathematicians became interested in gauge and string theories and from the physical side, physicists started to produce results that guided certain mathematical constructions. It can be said that a \textbf{new geometry} was born from those works. It is also important to mention that at that time Grothendieck's program of new foundations for algebraic geometry was in its highest pick, which also would have strong consequences in the new marriage between mathematics and physics although at the 70's that was not completely clear. Therefore, it is not exaggerated to say that, in a global sense, a new geometry was born during the 70's. Part of this new geometry has led to the development of the so-called Physical Mathematics.\\

In \cite{moore1, moore2} Moore addressed the question \emph{What is Physical Mathematics?}. The answer he provides, in his own words, goes as follows: \\

\emph{``The use of the term ``Physical Mathematics'' in contrast to the more traditional ``Mathematical Physics'' by myself and others is not meant to detract from the venerable subject of Mathematical Physics but rather to delineate a smaller subfield characterized by questions and goals that are often motivated, on the physics side, by quantum gravity, string theory, and supersymmetry, (and more recently by the notion of topological phases in condensed matter physics), and, on the mathematics side, often involve deep relations to
infinite-dimensional Lie algebras (and groups), topology, geometry, and even analytic number theory, in addition to the more traditional relations of physics to algebra, group theory,
and analysis. To repeat, one of the guiding principles is the goal of understanding the ultimate foundations of physics. Following the lessons of history, as so beautifully stressed by
Dirac, we may reasonably expect this to lead to important new insights in mathematics''.}\\

I am going to look at Moore's characterization in more detail, in order to extract what I consider are the main features of Physical Mathematics in his approach. \\

The first thing one should notice is the inversion of the terms mathematics and physics from the ``traditional Mathematical Physics'', pointing out how physics might be helpful in the articulation of mathematical ideas. In fact, Moore also adds (loc.cit.) that\\

\emph{``If a physical insight leads to a significant new result in mathematics,
that is considered a success. It is a success just as profound and notable as an experimental confirmation from a laboratory of a theoretical prediction of a peak or trough. For
example, the discovery of a new and powerful invariant of four-dimensional manifolds is a
vindication just as satisfying as the discovery of a new particle''.} \\

Changing the order of the terms mathematics and physics, i.e. \textbf{inverting the view}, is highly non-trivial. It is not a mere linguistic game, it brings deep conceptual implications for both mathematics and physics. In this approach physics plays a new role. Now, physics not only tries to understand the \emph{sensible reality} but also helps to understand the \emph{ mathematical reality} and this should count as part of the success of a physical theory. 
This claim has profound philosophical implications that I am going to start addressing at the end of this text. Of course, this is a long term project and I will only provide some glimpses here. I just want to mention that this could address some points of the criticism on string theory as a physical theory. From the point of view of Physical Mathematics, string theory is a genuine physical theory because it has enhanced our understanding of the mathematical world. It is clear that in order to accept this point of view, one really needs to think seriously on what the inversion process means, so one needs a \textbf{philosophy of inversion}. Otherwise, from the traditional point of view this can be considered just nonsense. \\

The second point, closely related to the inversion of the traditional approach of the relations between mathematics and physics, has to do with the fact that Physical Mathematics is a new field of knowledge. As we tried to argue above, Physical Mathematics emerged in the context of a new geometry. This is pretty clear in Moore's description. He mentions very sophisticated areas of mathematics and physics involved in the development of Physical Mathematics which are at the very heart of the new geometry I have mentioned. For a detailed account on the mathematics and physics involved in the works on Physical Mathematics and also for a complete list of references of relevant works see \cite{mooreetal1, mooreetal2, moore1}. However, Moore is a bit cautious and located Physical Mathematics as a subfield of Mathematical Physics. In this point I have to disagree with his description, since this is contradictory with the new vistas that the inversion of the terms mathematics and physics in the context of Physical Mathematics bring. \\

It is obvious that Physical Mathematics is closely related to Mathematical Physics, but it should not be considered a part of it, but rather an independent field with its own features. Both have a nonempty intersection, that is clear, but pointing out the differences is a crucial point. If the differences are not discussed in depth, I do not see the point of talking about Physical Mathematics. The emergence of Physical Mathematics \textbf{splits} and \textbf{augments} our knowledge and these two processes should be captured through the differentiation between Physical Mathematics and Mathematical Physics which is hidden in Moore's definition. \\

The third aspect I want to remark is the relation between Physical Mathematics and the foundations of physics. Here foundations should be understood in a general and ample sense, not only as the search of some definitive principles, nor as the search of a theory of everything. The foundations one should look for must be some mobile and elastic principles, not rigid and immutable ones. There is certain dynamism in the foundations of physics, and also in the foundations of mathematics, that must be captured by certain mobility and plasticity of the physical and mathematical entities. One needs to understand the emergency of virtualities both in mathematics and physics. Therefore, a \textbf{philosophy of virtuality} is needed in order to understand how new fundamental objects (physical and mathematical) emerge from the Physical Mathematical theories. \\

Let us summarize the discussion. We have observed the emergence of a new field of knowledge that invites rethinking the traditional relations between mathematics and physics. This field has been called Physical Mathematics, and motivated by the description proposed by one of its leading figures, I have argued that in order to clarify and understand better such description, and not only the description but also the power of this new field, it is necessary to introduce a philosophical approach that takes care of, at least, three process, namely, inverting, differentiating - subdivided in turn in two subprocesses, splitting and augmenting-, and the emergency of virtualities. Such a philosophical approach is already available and it constitutes part of the work of Gilles Ch\^atelet.

\section{Physical mathematics and Ch\^{a}telet's philosophy}

From his remarkable 1993 treatise \emph{Les enjeux du mobile. Math\'ematique, physique, philosophie}\footnote{I will use the English translation \emph{Figuring space. Mathematics, physics, philosophy} by Robert Shore and Muriel Zagha.} \cite{chatelet1}, it is pretty clear that Ch\^atelet's program concerns exploring the relations among mathematics, physics and philosophy as it is explicit in the title. Certainly, Ch\^atelet is not the first one who was interested in those relationships. One can trace works along those directions even in Plato, passing notably through Leibniz, Riemann, Poincar\'e or Lautman, just to mention some outstanding names. However, it is Ch\^atelet's point of view the closest one, as far as I can see, to the developments of Physical Mathematics, even though Ch\^atelet himself was not completely aware of the existence of this field.  \\

In his work Ch\^atelet uses the French expression \textbf{physico-math\'ematique}, which was translated into English as \textbf{physico-mathematics}. In his words: \\

\emph{``At first sight this physico-mathematics does appear to be an axiomatic giving precise form to the system of equivalence between mathematical concepts and physical concepts. To understand the revolutionary coup that installs this axiomatic is to discover the proximity of two horizons made up of virtual determinations which exceed the current set of explicit determinations and which still remain available for examination''.} \\

Reading this paragraph in 2025 it is inevitable to think of Physical Mathematics, as described by Moore. This actually validates, from my point of view, the translation as physico-mathematics and not as mathematical physics. Physico-mathematics is not Mathematical Physics, it is something new and Ch\^atelet was aware of that. The remarkable thing is that he wrote it in 1993! Hence, in that time, he already saw the new dynamics in the relations between mathematics and physics and envisioned the new field of Physical Mathematics, which he encoded in his physico-mathematics\footnote{In \cite{zaslow}, Zaslow uses the word \emph{physmatics} to describe the new links between mathematics and physics. According to him this word emphasizes that the links are profound and inseparable. However, in my point of view, the expression Physical Mathematics captures better those new links via the inversion of the terms in the expression Mathematical Physics.}. The role of the virtual emergencies is also very clear from his approach. It is in the emergence of the virtualities where the new relations between mathematics and physics take place and go beyond the traditional approach encapsulated in the idea of Mathematical Physics. It is also important to note that Ch\^atelet seems to see the necessity of inverting the term physics and mathematics in the expression Mathematical Physics. Along this lines he wrote:\\

\emph{``What gestures are involved in this `functioning', which establishes a higher form of continuity between `structure' and `reality'? Is it simply a case of `applying' mathematics to physics or rather of awakening the physical in mathematics?''} \\

\textbf{Awakening the physical in mathematics} is precisely one of the main goals of Physical Mathematics. I would like to add one more quote from \emph{Les enjeux du mobile} to make Ch\^atelet's call of a \textbf{new Mathematical Physics} (a.k.a., Physical Mathematics) clearer and also to point out the role of virtualities in the consolidation of this new field: \\

\emph{``Paradoxical result: it is the motion itself of the amplifying abstraction of mathematics that governs their incarnation as physical beings: the more `abstract' mathematics is, the better it works in application. Thus, to establish a new mathematical physics is to be capable of recognizing to what extent such an effort of radical autonomy in mathematics necessarily involves the horizon of virtualities of physics''.}
\\
 
The recognition of a new Mathematical Physics also has the consequence of establishing precise differences between the traditional Mathematical Physics and the new one. Of course, this requires understanding more sophisticated mathematics and physics. Therefore, Ch\^atelet focused on gauge theories to work out his program. In the posthumous \emph{L'enchantement du virtuel} \cite{chatelet2}, in particular, in the essays \emph{Le potentiel d\'emoniaque. Aspects philosophiques et physiques de la th\'eorie de jauge} and \emph{La physique mate\'ematique comme projet. Un exemple : la grande unification des forces}, Ch\^atelet takes modern and sophisticated developments in physics to undertake his philosophical analysis. In fact, gauge theories lie at the heart of Physical Mathematics\footnote{See \cite{atiyah} for a nice layman introduction to gauge theories and its relation with modern physics and mathematics.}, as same as the problem of the unification of forces under the disguise of string theory. The following quote from \emph{La physique math\'ematique comme projet. Un exemple : la grande unification des forces} is very illustrative: \\

\emph{``Cest probablement comme nouveau mode d'instauration du physico-math\'ematique lui-m\'eme que la th\'eorie de jauge est r\'evolutionnaire. Toute physique math\'ematique est fond\'ee par un protocole d'articulation du sens g\'eom\'etrique et du sens physique''.} \\

And again the emergence of the virtuality is explicit, but now the reflection focuses on the gauge theories. Here Ch\^atelet's vision is remarkable. He is putting the finger precisely on one of the central theories in Physical Mathematics. Of course, by the time Ch\^atelet wrote his essay, gauge theories were already a well established theory in the mathematical and physical worlds, but what is interesting to pointing out is how Ch\^atelet as ``philosopher'' sees the need to discuss the current developments of the science, which is not the most common attitude in the philosophical community. Using gauge theory for a philosophical analysis, in relation with virtualities, Ch\^atelet wrote : \\

\emph{``Si maintenant la rencontre du physique et du mathématique est comprise explicitement comme l'acte de traduction qui identifie un rassemblement concret de fibres par transport parallèle et un potentiel-vecteur créateur de particules virtuelles, véhicules de l'interaction comme e'est le cas pour la théorie de jauge, le malaise s'évanouit''.} \\

In \cite{zalameachat}, Zalamea analyses an unpublished worksheet of Ch\^atelet and detects three important processes in his methodology, which are more or less clear in the published work, but complete evident in the worksheet studied by Zalamea. The mentioned processes are inverting, splitting and augmenting. In particular, Zalamea points out inverting, splitting and augmenting are the \emph{``three key gestures to understand the mathematical fabric''}\footnote{Along the text I have been referring of inverting, splitting and augmenting as processes, but in a more consistent and precise way with Ch\^atelet's philosophy, they must be called \textbf{gestures}, as Zalamea rightly does. In fact, in one the Ch\^atelet's quotes above, one can see that he is asking for what gestures are needed for understanding the physical in mathematics. On the other hand, I also fused splitting and augmenting in one gesture I call differentiating. With this fusion I want to emphasize that in order to make differences between fields of knowledge we first split the knowledge and then it naturally augments.}. In fact, according to the study presented in this text, they are also the three key gestures to understand Physical Mathematics. With these three gestures one can build the philosophy of inversion, differentiation and virtuality that was required to analyse Moore's description of Physical Mathematics. This is quite natural since, as it has been discussed, Ch\^atelet's program from the very beginning was advocating for the construction of new relations between mathematics and physics of the sort presented in what is now called Physical Mathematics.

\section{Physical Mathematics as a project. An example: mirror symmetry} 

The title of this section evokes Ch\^atelet's essay \emph{La physique math\'ematique comme projet. Un exemple : la grande unification des forces} \footnote{It is curious how Ch\^atele uses here the more traditional \emph{physique math\'ematique} rather than his \emph{physico-math\'ematique}. }. However, I do not pretend to mimic Ch\^atelet's impressive work, I just want to discuss in a bird's-eye view  a future direction in the  study of the philosophy of Physical Mathematics. The few line that come are only an invitation to a vast program, the main study case for taking Physical Mathematics as a project. \\

Around the time Ch\^atelet was writing his \emph{Les enjeux du mobile}, a remarkable breakthrough along the lines of using physical insights to produce new mathematics was happening. Candelas, de la Ossa, Green and Parkes in 1991 wrote a remarkable paper about how to use mirror symmetry ideas (from a physical point of view) to tackle old and interesting problems in enumerative geometry, see \cite{candelas}. Of course, the paper goes beyond this and opened fruitful developments both in mathematics and physics, that are still taking place and are not yet fully understood. There are only two years of difference between this work and Ch\^atelet's book, short time for Ch\^atelet to be aware of this, in particular, because he was not an expert in the field\footnote{In addition, unfortunately, he passed away just after few years of the publication of \emph{Les enjeux du mobile}, so we will never be able to know what kind of things Ch\^atelet might have done knowing about the power of mirror symmetry and its influence in mathematics and physics.}. Nevertheless, what it is interesting for me is this neighbourhood of two years\footnote{Also in 1994, in his ICM talk, Maxim Kontsevich launched his own program for understanding mirror symmetry from a categorical point of view, see \cite{kontsevich}. This program is known as homological mirror symmetry and it has been of great influence in both mathematics and physics. Nonetheless, it is begging for a philosophical treatment. A first attempt of this can be found in \cite{zalamea}.} where many important things from  mathematics, physics and philosophy happened in such a way that changed our perception of how these fields should be related. In fact, from my point of view, the paper by Candelas, de la Ossa, Green and Parkes is a masterpiece of Physical Mathematics, and the whole mirror symmetry should take now the place that gauge theories had in Ch\^atelet's work. As Atiyah, Dijkgraaf and Hitchin \cite{atiyah2} wrote:\\

\emph{``Mirror symmetry is a good example of a more fundamental influence of physics on geometry, one which involves a significant change of viewpoint on the part of the pure mathematician.''}\\

Mirror symmetry is still an intriguing phenomenon. Arising in physics as a duality in conformal field theories, it can be seen in mathematics as a duality between symplectic and complex geometry\footnote{This is two dimensional mirror symmetry which in its original formulation is a duality between a pair of Calabi-Yau manifolds interchanging their symplectic and complex geometry. This phenomenon can be extended beyond Calabi-Yau's and in fact, we have, for example, Fano-Landau-Ginzburg model mirror correspondence too. On the other hand, one can consider not a two dimensional theory but a three dimensional theory and we have the so called three dimensional mirror symmetry, which mathematical incarnation is a symplectic duality between certain holomorphic symplectic manifolds. In addition, the geometric Langlands duality can be understood as a kind of four dimensional mirror symmetry. This different aspects of mirror symmetry will be discussed in detail in \cite{cruz1}. In this text, I will consider mirror symmetry as two dimensional mirror symmetry.}. To be close to the physics origin of mirror symmetry, I will discuss mirror symmetry for Calabi-Yau 3-folds, but the restriction of the dimension can be removed. Given a Calabi-Yau manifold $X$, there are two natural geometrical structures one can attach to it, namely a symplectic structure and a complex structure. Mirror symmetry claims that for a Calabi-Yau 3-fold $X$, there exists another Calabi-Yau 3-fold $X^\vee$ such that the symplectic geometry of $X$ is related to the complex geometry of $X^\vee$ and vice-versa. In fact, mirror symmetry should be formulated in terms of families of Calabi-Yau 3-folds but this is a technical point which is not important for our discussion. The symplectic side is called the A-model side and the complex side is the B-model side. \\

What Candelas, de la Ossa, Green and Parkes did was to take the problem of counting rational curves for a specific Calabi-Yau 3-fold $X$ (the quintic hypersurface inside $\mathbb{CP}^4$), a problem attached with the geometry of the A-model, and working with the geometry of the B-model, i.e. the complex geometry of the mirror (family) of $X$, they found the solution for the counting problem, thanks to the duality given by the mirror phenomenon. It is interesting to see how going from the A-model to the B-model or the other way around, might be seen as an inverting gesture since one is inverting the side of the problem. Why this works is still a mystery but perhaps a conceptual approach using the gesture of inverting could clarify something, at least from a philosophical point of view. In any case, mirror symmetry has proven to be a very powerful theory in both mathematics and physics and it goes beyond the application to enumerative problems.\\

In mirror symmetry the emergence of virtualities is pretty clear\footnote{Explicit virtualities, like the \emph{virtual fundamental class} and \emph{virtual dimensions} in Gromov-Witten theory (a study of certain geometric invariants arising in the A-model), arise quite naturally and they are needed for the right mathematical foundations of some of the objects involved in the mirror phenomenon.}, so taking care of such emergencies from a philosophical perspective could be of benefit for understanding mirror symmetry. As Zalamea points out in \cite{zalameachat}: \emph{``ramified virtualities may help to understand the `blind points' of creativity''}. The lack of understanding on why the mirror phenomenon occurs seems to be one of these blind points. \\

From a purely mathematical point of view mirror symmetry also has the advantage of integrating different branches of mathematics. As we already mentioned, symplectic and complex geometry come at first sight but also integrable systems, derived categories, geometric representation theory, modular forms, L-functions, ordinary differential equations and isomonodromic problems and noncommutative geometry arise naturally in considering the mathematics around mirror symmetry. However, from a philosophical analysis point of view, and also thinking in potential applications in physics, there is one area of mathematics that it might be specially relevant to understand the mirror dualities but that is not part of the mainstream of the community working in mirror symmetry, not at least up to now. That branch is model theory. \\

Recently model theory has called the interest of many researches because of its applications, under the disguise of o-minimality, in Hodge theory and arithmetic geometry. In a certain sense, this is also a manifestation of the importance of the tame geometry program of Grothendieck, since o-minimality is generally considered a way to approach tame geometry. Due to this proximity, specially with Hodge theory, o-minimality looks relevant to understand certain aspects of mirror symmetry, for instance periods maps. However, I do not want to call the attention to o-minimality but rather to a different approach of model theory in mirror symmetry via categoricity and its variations. \\

Zilber's program, see \cite{zilber1,zilber2}, is a proposal of using the idea of logically perfect structures, in order to tackle problems in geometry, in particular, arithmetic geometry and noncommutative geometry, and also to study the foundations of physics. According to the discussion presented in this text, one can put Zilber's program inside the Physical Mathematical enterprise. It is not completely clear how Zilber's program fits inside mirror symmetry, but there are some intuitions on how to proceed. Here, I just want to mention that this project is, in fact, part of a bigger one aiming for looking further in the physical mathematics around mirror symmetry. In fact, being a bit optimistic, mirror dualities could be a manifestation of the more general syntactic/semantics dualities appearing in Zilber's approach, in the same way as one can see mirror symmetry as mediating in the polarity between the discrete and the continuum\footnote{In a more general picture, mirror symmetry can be seen as a \emph{partially solving} some of the important polarities in mathematics like discrete/continuum, local/global, one/multiple, structure/deformations.}. \\

Certainly undertaking the project of analysing mirror symmetry as a part of Physical Mathematics, with detailed technical discussions from both mathematics and physics and with the technical philosophical framework provided by Ch\^atelet's philosophy, looks like an important project. Also, how to use the insights coming from philosophy and physics to produce new mathematics seems to be a good challenge. This will show the relevance of Ch\^atelet's ideas and would be a partial realization of Grothendieck's vision where the fundamental questions about the structure of the universe need a mathematical heart, physical insights and philosophical openness\footnote{It could be the case that this vision is also implicit in the work of many other important thinkers (mathematicians, philosophers, physicist), but make it explicit is important since it might help to reveal the mysteries of the activity of creation.}. \\

\section*{Acknowledgements}

I would like to express my great gratitude to Fernando Zalamea who encourages me to write this paper. He introduced me to Gilles Ch\^atelet's philosophy during his long term \emph{Seminario de Filosof\'ia Matem\'atica} (Philosophy of Mathematics Seminar) at Universidad Nacional de Colombia and also sent me a copy of \cite{zalameachat}, which has been essential for the analysis presented in this text.\\
I would also like to thank Michael Harris, Juliette Kennedy, Alejandra Mu\~noz, Motohico Mulase, Andr\'es Villaveces and Boris Zilber for stimulating conversations about mathematics, physics and philosophy. Most of our conversations helped me to give shape to this essay.\\
Part of this paper was written during the time I was a CAS-PIFI fellow as visiting professor of the University of Science and Technology of China. I thank the Chinese Academy of Science and the University of Science and Technology of China for the financial support and excellent working conditions.




\end{document}